\documentclass[11pt]{amsart}
\usepackage[latin1]{inputenc}
\usepackage{a4}
\usepackage{amssymb}
\usepackage{amsmath}
\begin{document}
\setlength{\parindent}{1.2em}
\def\COMMENT#1{}
\def\TASK#1{}
\def\noproof{{\unskip\nobreak\hfill\penalty50\hskip2em\hbox{}\nobreak\hfill%
       $\square$\parfillskip=0pt\finalhyphendemerits=0\par}\goodbreak}
\def\endproof{\noproof\bigskip}
\newdimen\margin   
\def\textno#1&#2\par{%
   \margin=\hsize
   \advance\margin by -4\parindent
          \setbox1=\hbox{\sl#1}%
   \ifdim\wd1 < \margin
      $$\box1\eqno#2$$%
   \else
      \bigbreak
      \hbox to \hsize{\indent$\vcenter{\advance\hsize by -3\parindent
      \sl\noindent#1}\hfil#2$}%
      \bigbreak
   \fi}
\def\proof{\removelastskip\penalty55\medskip\noindent{\bf Proof. }}
\def\enddiscard{}
\long\def\discard#1\enddiscard{}
\newtheorem{firstthm}{Proposition}
\newtheorem{thm}[firstthm]{Theorem}
\newtheorem{prop}[firstthm]{Proposition}
\newtheorem{lemma}[firstthm]{Lemma}
\newtheorem{cor}[firstthm]{Corollary}
\newtheorem{problem}[firstthm]{Problem}
\newtheorem{defin}[firstthm]{Definition}
\newtheorem{conj}[firstthm]{Conjecture}
\newtheorem{theorem}[firstthm]{Theorem}
\newtheorem{claim}[firstthm]{Claim}
\def\eps{{\varepsilon}}
\def\N{\mathbb{N}}
\def\R{\mathbb{R}}
\def\K{\mathcal{K}}
\def\vH{\vec{H}}
\def\vG{\vec{G}}

\title{A note on complete subdivisions in digraphs of large outdegree}
\author{Daniela K\"uhn \and Deryk Osthus \and Andrew Young}
\date{}
\maketitle \vspace{-.8cm}
\begin{abstract} \noindent
Mader conjectured that for all $\ell$ there is an integer  $\delta^+(\ell)$ such that
every digraph of minimum outdegree at least $\delta^+(\ell)$ contains a subdivision of a
transitive tournament of order $\ell$. In this note we observe that if the minimum outdegree of
a digraph is sufficiently large compared to its order
then one can even guarantee a subdivision of a large complete digraph.
More precisely, let $\vG$ be a digraph of order~$n$ whose minimum outdegree is at least $d$. 
Then $\vG$ contains a subdivision of a complete
digraph of order $\lfloor d^2/(8n^{3/2}) \rfloor$. 
\end{abstract}

\section{Introduction}\label{intro}

A fundamental result of Mader~\cite{Mader1967} states that
for every integer $\ell$ there is a smallest $d=d(\ell)$ so that
every graph of average degree at least
$d$ contains a subdivision of a complete graph on $\ell$ vertices.
Bollob\'as and Thomason~\cite{BT1998PCMEHTCS} as well as
Koml\'os and Szemer\'edi~\cite{KS1996TCGII}
showed that $d(\ell)$ is quadratic in $\ell$.
In~\cite{M1985DLCD}, Mader made the following conjecture,
which would provide a digraph analogue of these results
(a transitive tournament is a complete graph whose edges are oriented 
transitively).
\begin{conj}[\textbf{Mader~\cite{M1985DLCD}}]\label{Maderconj}
For every integer $\ell>0$ there is a smallest integer $\delta^+(\ell)$ 
such that every digraph~$\vG$ with minimum outdegree at least  
$\delta^+(\ell)$ contains a subdivision of the transitive tournament on~$\ell$
vertices.
\end{conj}
It is easy to see that $\delta^+(\ell)=\ell-1$ for $\ell \leq 3$.
Mader~\cite{M1996OTTOFDOT} showed that $\delta^+(4)=3$. 
Even the existence of $\delta^+(5)$ is not known.
One might be tempted to conjecture that large minimum outdegree
would even force the existence of a subdivision of a large complete digraph
(a complete digraph has a directed edge from $v$ to $w$ for any
ordered pair $v,w$ of vertices). However, for all~$n$ 
Thomassen~\cite{T1985ECDG} constructed a digraph on $n$ vertices
whose minimum outdegree is at least $\frac{1}{2}\log_2 n$ 
but which does not contain an even directed cycle
(and thus no complete digraph on $3$ vertices).
The additional assumption of large minimum indegree in Conjecture~\ref{Maderconj}
does not help either. Mader~\cite{M1985DLCD} modified the construction in~\cite{T1985ECDG}
to obtain digraphs having arbitrarily large minimum indegree
and outdegree without a subdivision of a complete digraph on $3$ vertices.

The fact that one certainly cannot replace the minimum outdegree in Conjecture~\ref{Maderconj} 
by the average degree is easy to see:
consider the complete bipartite graph with equal size vertex classes and orient all edges from the
first to the second class. The resulting digraph $\vec{B}$
has average degree $|\vec{B}|/2$ but not even a directed cycle
or a transitive tournament on $3$ vertices.
(On the other hand, 
Jagger~\cite{J1998EDRTCS} showed that if the average degree of a
digraph~$\vec{G}$ is a little larger than
$|\vec{G}|/2$, then $\vec{G}$ does contain a subdivision of a large complete digraph.)

So in some sense, the above examples and constructions show that
Conjecture~\ref{Maderconj} is the only possible analogue of the result
in~\cite{Mader1967} mentioned above. Our main result is that if the minimum
outdegree of a digraph is sufficiently large compared to its order,
then Conjecture~\ref{Maderconj} is true. In fact, we show that in this case,
one can even guarantee a subdivision of a complete digraph.

\begin{thm}\label{thm:digraph}
Let $\vG$ be a digraph of order~$n$ whose minimum outdegree is at least~$d$.
Then~$\vG$ contains a subdivision of the complete
digraph of order $\lfloor d^2/(8n^{3/2}) \rfloor$.
\end{thm}
Note that the bound is nontrivial as soon as $d$ is a little larger than $n^{3/4}$.
Also, recall that the result of Thomassen~\cite{T1985ECDG} mentioned above implies that
we cannot have a subdivision of a complete digraph of order at least~$3$
if $d\le \frac{1}{2}\log_2 n$. Furthermore, note that if $d=cn$, then
Theorem~\ref{thm:digraph} guarantees a subdivision of a complete
digraph of order $\lfloor c' \sqrt{n} \rfloor$, where $c'=c^2/8$.
It is easy to see that this is best possible up to the value of $c'$ 
(consider the complete bipartite digraph with vertex classes of equal size).

The main ingredient in the proof of Theorem~\ref{thm:digraph} is
Lemma~\ref{lemma:subgraph}. It states that if~$\vG$ has~$n$ vertices
and its minimum outdegree is $\gg \sqrt{n}$, then
$\vG$ has a subdigraph $\vH$ which is highly connected in the
following sense: if $x$ is any vertex of~$\vH$ and $y$ is a vertex of large
indegree, then there are many internally
disjoint dipaths from $x$ to $y$ in~$\vH$.
Lemma~\ref{lemma:subgraph} also guarantees the existence of
many such vertices~$y$. For undirected graphs, there is a much stronger
result of Mader~\cite{Mader72} which implies that every graph of minimum
degree at least $4k$ has a $k$-connected subgraph.
Since a digraph version of this result is not known,
Lemma~\ref{lemma:subgraph} may be of independent interest.
There are also several related results of Mader~\cite{M1985DLCD,M1995EVLCDLO} 
which investigate the existence of pairs of vertices
with large local connectivity in digraphs of large minimum outdegree.
The proof of Lemma~\ref{lemma:subgraph} is quite elementary:
if the current subdigraph~$\vec{H}$ does not satisfy the requirements,
then we can use Menger's theorem to find a significantly smaller subdigraph whose minimum 
outdegree is almost as large as that of~$\vec{H}$. Since this means that the density of the 
successive subdigraphs increases, this process must eventually terminate.

\section{Proof of Theorem~\ref{thm:digraph}}

Before we start with the proof of Theorem~\ref{thm:digraph} let us introduce
some notation.
The digraphs~$\vG$ considered in this note do not contain loops and between
any ordered vertex pair $x,y\in \vG$ there is at most one edge from~$x$ to~$y$.
(There might also be another edge from~$y$ to~$x$.)
We denote by $\delta^+(\vG)$ the minimum outdegree of a digraph~$\vG$
and by~$|\vG|$ its order.
We write $d^+_{\vG}(x)$ for the outdegree of a vertex~$x\in\vG$ and
$d^-_{\vG}(x)$ for its indegree. A digraph~$\vH$ is a \emph{subdivision} of~$\vG$
if~$\vH$ can be obtained from~$\vG$ by replacing each edge~$\vec{xy}\in\vG$ with
a dipath from~$x$ to~$y$ such that all these dipaths are internally disjoint for
distinct edges. The vertices of~$\vH$ corresponding to the vertices of~$\vG$ are
called \emph{branch vertices}.
 
Given two vertices~$x$ and~$y$ of a digraph~$\vG$, we define $\kappa_{\vG}(x,y)$
to be the largest integer $1\le k\le |\vG|-2$ such that $\vG-S$ contains
a dipath from~$x$ to~$y$ for every vertex set $S\subseteq V(\vG)\setminus\{x,y\}$
of size~$<k$. We define $\kappa_{\vG}(x,y):=0$ if $\vG$ does not contain
a dipath from~$x$ to~$y$. We will use the following version of Menger's theorem
for digraphs.

\begin{thm}[\textbf{Menger's theorem for digraphs}]\label{thm:diMenger} 
Let $x$ and $y$ be vertices of a digraph~$\vG$ such that $\kappa_{\vG}(x,y)\ge k$.
Then~$\vG$ contains~$k$ internally disjoint dipaths from~$x$ to~$y$.
\end{thm}

As mentioned above, the main step in the proof of Theorem~\ref{thm:digraph} is to find
a subdigraph~$\vH$ of~$\vG$ such that the minimum outdegree of~$\vH$ is still
large and such that every vertex of~$\vH$ sends many internally disjoint dipaths to
each vertex of~$\vH$ which has large indegree. 

\begin{lemma}\label{lemma:subgraph}
Every digraph~$\vG$ of order~$n$ with $\delta^+(\vG)\ge d$ contains
a subdigraph~$\vH$ such that
\begin{itemize}
\item[{\rm (i)}] $\delta^+(\vH)> d/2$,
\item[{\rm (ii)}] $\kappa_{\vH}(x,y)\ge d^2/(4n)$ for all pairs
$x,y\in V(\vH)$ with $d^-_{\vH}(y)\ge d/2$,
\item[{\rm (iii)}]  at least $d^2/(4n)$
vertices of~$\vH$ have indegree at least~$d/2$ in~$\vH$.
\end{itemize}
\end{lemma}
\proof
Put $$\alpha:=\frac{d}{n} \ \ \ \ \text{and}\ \ \ \
\alpha':=\frac{d^2}{4n^2}=\frac{\alpha^2}{4}.$$
By Theorem~\ref{thm:diMenger} we may assume that $\kappa_{\vG}(x,y)<\alpha'n$
for some vertices $x,y$ of~$\vG$ with $d^-_{\vG}(y)\ge d/2$.
Otherwise we could take $\vH:=\vG$. (It is easy to check that~$\vH$ then
also satisfies condition~(iii) of the lemma.)%
     \COMMENT{Let $\ell$ denote the number of vs of degree $\ge d/2$.
We have $dn\le \delta^+(\vG)n\le \ell n+(n-\ell)d/2$ and so
$dn/2\le \ell(n-d/2)\le \ell n$ ie $\ell\ge d/2\ge d^2/(4n)$.} 
Let $S\subseteq V(\vG)\setminus\{x,y\}$ be a set of size $<\alpha'n$ such that~$\vG-S$ does
not contain a dipath from~$x$ to~$y$. Let $Y$ be the set of all those vertices~$z$
for which $\vG-S$ contains a dipath from~$z$ to~$y$. Then~$Y\cup S$ contains~$y$
as well as all the at least~$d/2=\alpha n/2$ inneighbours of~$y$. Let~$C$ denote the component
of the undirected graph corresponding to $\vG-(Y\cup S)$ which contains~$x$.
Let~$\vG_1$ be the subdigraph of~$\vG$ induced by all vertices in~$C$.
Then $|\vG_1|\le n-|Y\cup S|< (1-\alpha/2)n$. Moreover, note that there exists
no edge directed from a vertex of~$\vG_1$ to a vertex outside $V(\vG_1)\cup S$.
Thus
\begin{equation}\label{eqminoutdeg}
\delta^+(\vG_1)\ge \delta^+(\vG)-|S|> (\alpha-\alpha')n.
\end{equation}
If~$\vG_1$ does not satisfy condition~(ii) of the lemma we again apply Theorem~\ref{thm:diMenger}
to obtain a subdigraph $\vG_2\subseteq \vG_1$. We continue in this fashion
until we obtain a subdigraph~$\vG_r$ which satisfies condition~(ii).
We will show that~$\vG_r$ also satisfies~(i) and~(iii). Put $\vG_0:=\vG$,
$$\delta_i:=\frac{\delta^+(\vG_i)}{|\vG_i|}\ \ \ \ \text{and}\ \ \ \
\gamma_{i-1}:=\frac{|\vG_{i-1}|}{|\vG_i|}$$
for all $i\le r$.
Similarly as in~(\ref{eqminoutdeg}) it follows that
\begin{equation}\label{eqminout1}
\delta^+(\vG_i)=\delta_i|\vG_i|\ge \delta_{i-1}|\vG_{i-1}|-\alpha'n\ge (\alpha-i\alpha')n.
\end{equation}
Thus
$\delta_i\ge \delta_{i-1}\gamma_{i-1}-\alpha' n/|\vG_i|=
\delta_{i-1}\gamma_{i-1}-\alpha'\prod_{j=0}^{i-1} \gamma_j$.
Using this inequality and induction on~$i$ one can show that%
     \COMMENT{Indeed, it is clear that \(\delta_1\geq
\gamma_0\delta_{0}-\gamma_0\alpha' =
\gamma_0(\alpha-\alpha')\). Moreover
$$
\delta_{l+1} \geq
\gamma_{l}\delta_{l}-\left(\prod_{j=0}^{l}\gamma_{j}\right)\alpha'
\ge \gamma_{l}\left(\left(\prod_{j=0}^{l-1}\gamma_{j}\right)(\alpha-l\alpha')\right)-
\left(\prod_{j=0}^{l}\gamma_{j}\right)\alpha'
=\left(\prod_{j=0}^{l}\gamma_{j}\right)(\alpha-(l+1)\alpha')
$$}
\begin{equation}\label{eqminout2}
\delta_i\ge (\alpha-i\alpha')\prod_{j=0}^{i-1} \gamma_j=(\alpha-i\alpha')\frac{n}{|\vG_i|}.
\end{equation}
Since we delete at least $d/2=\alpha n/2$ vertices when going from
$\vG_{i-1}$ to $\vG_i$ (namely the inneighbours
of the vertex playing the role of~$y$), we have that $|\vG_r|\le n-r\alpha n/2$.
In particular this shows that $r<2/\alpha$. However, since~(\ref{eqminout2}) implies that
$1>\delta_r\ge (\alpha-r\alpha')/(1-r\alpha/2)$ we even have%
     \COMMENT{If $r<2/\alpha$ then the inequality is equivalent to
$1-r\alpha/2\ge \alpha-r\alpha'$ ie $1-\alpha\ge r(\alpha/2-\alpha')$.
Note that the 2nd bound is indeed stronger: $2/\alpha\ge (1-\alpha)/(\alpha/2-\alpha')=
(1-\alpha)/(\alpha/2-\alpha^2/4)$ since this is equivalent to
$1-\alpha/2\ge 1-\alpha$.}
$r<(1-\alpha)/(\alpha/2-\alpha')$. Thus%
    \COMMENT{$\alpha-\frac{1-\alpha}{\alpha/2-\alpha'}\alpha'=
\alpha-\frac{1-\alpha}{2/\alpha-1} =\alpha-\frac{\alpha(1-\alpha)}{2-\alpha}=
\frac{2\alpha-\alpha^2-\alpha+\alpha^2}{2-\alpha}=
\frac{\alpha}{2-\alpha}$.}
\begin{equation}\label{eqminout3}
\delta^+(\vG_i)\stackrel{(\ref{eqminout1})}{\ge}
(\alpha -r \alpha')n\ge \left(\alpha-\frac{1-\alpha}{2/\alpha-1}\right)n
=\frac{\alpha n}{2-\alpha}>\frac{d}{2}.
\end{equation}
Altogether this shows that~$\vG_r=:\vH$ satisfies conditions~(i) and~(ii) of the lemma.
To check that $\vH$ also satisfies condition~(iii)
let $\ell$ denote the number of vertices of indegree $\ge d/2$ in~$\vH$.
Then $$\frac{\alpha n|\vH|}{2-\alpha}\stackrel{(\ref{eqminout3})}{\le}
\delta^+(\vH)|\vH|\le |\vH|\frac{d}{2}+\ell |\vH|,
$$
which implies that $\ell\ge \alpha d/(4-2\alpha)\ge d^2/(4n)$, as%
     \COMMENT{Indeed, $\ell\ge \frac{\alpha n}{2-\alpha}-\frac{d}{2}=
d(1/(2-\alpha)-1/2)=d\frac{2-(2-\alpha)}{2(2-\alpha)}
=\frac{\alpha d}{4-2\alpha}.$}
required.
\endproof

\removelastskip\penalty55\medskip\noindent{\bf Proof of Theorem~\ref{thm:digraph}. }
Let $\ell:=\lfloor d^2/(8n^{3/2})\rfloor$.
We first apply Lemma~\ref{lemma:subgraph} to obtain a subdigraph
$\vH\subseteq \vG$ as described there. We pick a set~$X\subseteq V(\vH)$ of
$\ell$ vertices having indegree~$\ge d/2$ in~$\vH$. (Such a set~$X$ exists
by condition~(iii) of Lemma~\ref{lemma:subgraph}.) $X$ will be the set of
our branch vertices. For every pair $x,y\in X$ there exist
at least~$d^2/(4n)$ internally disjoint dipaths from~$x$ to~$y$. Thus the
average number of inner vertices on such a path is at most~$4n^2/d^2$.
Hence~$\vH$ contains at least~$d^2/(8n)$ internally disjoint dipaths from~$x$ to~$y$
such that each of these has at most~$8n^2/d^2$ inner vertices.
Let us call such a dipath~\emph{short}.
This shows that we can connect all pairs of branch vertices greedily (in both directions)
by choosing each time a short dipath which is internally disjoint from all the short dipaths
chosen before. In each step we destroy at most~$8n^2/d^2$ further dipaths.
But $(|X|^2-1)8n^2/d^2<8\ell^2n^2/d^2\le d^2/(8n)$, so we can connect
all pairs of branch vertices by short dipaths.
\endproof

\medskip

{\footnotesize \obeylines \parindent=0pt

Daniela K\"{u}hn, Deryk Osthus \& Andrew Young
School of Mathematics
University of Birmingham
Edgbaston
Birmingham
B15 2TT
UK
}

{\footnotesize \parindent=0pt

\it{E-mail addresses}:
\tt{\{kuehn,osthus,younga\}@maths.bham.ac.uk}}


\begin{thebibliography}{10}
\bibitem{BT1998PCMEHTCS} B.~Bollob\'{a}s and A.~Thomason, Proof of
a conjecture of Mader, Erd\H{o}s and Hajnal on topological
complete subgraphs, \emph{European Journal of
Combinatorics}~\textbf{19} (1998), 883--887.
\bibitem{J1998EDRTCS} C.~Jagger, Extremal Digraph Results for
Topological Complete Subgraphs, \emph{European Journal of
Combinatorics}~\textbf{19} (1998), 687--694.
\bibitem{KS1996TCGII} J.~Koml\'{o}s and E. Szemer\'{e}di,
Topological Cliques in Graphs II, \emph{Combinatorics, Probability
and Computing}~\textbf{5} (1996), 70--90.
\bibitem{Mader1967} W.~Mader, Homomorphieeigenschaften und mittlere Kantendichte
von Graphen, \emph{Math.~Annalen}~\textbf{174} (1967), 265--268.
\bibitem{Mader72} W.~Mader, Existenz $n$-fach zusammenh\"angender Teilgraphen
in Graphen gen\"ugend gro\ss er Kantendichte, \emph{Abh.~Math.~Sem.~Univ.~Hamburg}
\textbf{37}~(1972), 86--97.
\bibitem{M1985DLCD} W.~Mader, Degree and Local Connectivity
in Digraphs, \emph{Combinatorica}~\textbf{5} (1985), 161--165.
\bibitem{M1996OTTOFDOT} W.~Mader, On Topological Tournaments of
order 4 in Digraphs of Outdegree 3, \emph{Journal of Graph
Theory}~\textbf{21} (1996), 371--376.
\bibitem{M1995EVLCDLO} W.~Mader, Existence of vertices of local
connectivity \(k\) in digraphs of large outdegree,
\emph{Combinatorica}~\textbf{15} (1995), 533--539.
\bibitem{T1985ECDG} C.~Thomassen, Even Cycles in Directed Graphs,
\emph{European Journal of Combinatorics}~\textbf{6} (1985),
85--89.
\end{thebibliography}
\end{document}